\documentclass[11pt,leqno]{book}
\usepackage{amssymb,amsmath,amsthm,rotating}

\numberwithin{equation}{chapter}

\newcommand{\n}{\noindent}
\newcommand{\bb}[1]{\mathbb{#1}}
\newcommand{\cl}[1]{\mathcal{#1}}

\newcommand{\ovl}{\overline}
\newcommand{\sst}{\scriptstyle}

\theoremstyle{plain}
\newtheorem{lem}{Lemma}[chapter]

\newtheorem{thm}{Theorem}[chapter]

\newtheorem*{cor}{Corollary}

\theoremstyle{definition}
\newtheorem{defn}{Definition}[chapter]

\newtheorem{rem}{Remark}[chapter]

\begin{document}

\title{Operator Theory and Complex Geometry}

\author{Ronald G.~Douglas}

\date{}
\maketitle

\tableofcontents
\thispagestyle{empty}

\clearpage\thispagestyle{empty}\cleardoublepage

\chapter{FIRST LECTURE}\label{chp1}

\section{Introduction}\label{sec1.1}

\footnotetext{This is an informal writeup of a series of three lectures given at the Fourth Advanced Course in Operator Theory and Complex Analysis, Sevilla, 2007.\newline 
2000 Mathematics Subject Classification. 46E22, 46M20, 47B32, 32B99, 32L05.\newline 
Key words and phrases. Hilbert modules, \v Silov modules, kernel Hilbert spaces, invariant subspaces, isometries, holomorphic structure, localization.}

We discuss an approach to the study of bounded linear operators on a complex Hilbert space which involves concepts and techniques from complex geometry. Although the main goal is to develop an effective approach to the study of multivariate operator theory, the methods are also useful for the single variable case and we will use concrete examples from the latter to illustrate the theory.

We begin by recalling a basic result from linear algebra on finite dimensional Hilbert space, that the number of distinct eigenvalues of a matrix can't exceed the dimension of the space. It is not surprising that on an infinite dimensional Hilbert space an operator can have infinitely many distinct eigenvalues. However, the ``size'' or the cardinality of such a set is, perhaps, unexpected. In particular, consider the basic separable Hilbert space, $\ell^2({\bb Z}_+)$, of one-sided sequences $\{a_n\}_{n\in {\bb Z}_+}$, of complex numbers for ${\bb Z}_+$ the non-negative integer, which are square-summable. On $\ell^2({\bb Z}_+)$ a diagonal operator can have $\aleph_0$ distinct eigenvalues. However, a bounded operator on $\ell^2({\bb Z}_+)$ can have more.

Let $U_+$ be the unilateral shift on $\ell^2({\bb Z}_+)$; that is, $U_+\{(a_0,a_1,a_2,\ldots)\} = (0,a_0,a_1,\ldots)$. Then $U^*_+$, its adjoint the backward shift, has the property that $U^*_+k_\lambda = \bar\lambda k_\lambda$ for $\lambda$ in the open unit disk ${\bb D}$. Here $k_\lambda = (1,\bar\lambda, \bar\lambda^2,\ldots)$, which is square summable for $|\lambda|<1$. Thus $U^*_+$ has distinct eigenvalues corresponding to ${\bb D}$ which has the cardinality of the continuum. The unilateral shift and its adjoint were introduced by von Neumann  in his classic study \cite{Neu} of symmetric operators.

The family $\{k_\lambda\}$ is not only continuous in $\lambda$ but is an anti-holomorphic vector-valued function $\lambda\to k_\lambda$ on ${\bb D}$. (The introduction of the complex conjugate is standard so that there is holomorphicity for $U_+$ and anti-holomorphicity for $U^*_+$.) We will use this property to study the operator as follows. If we consider $\ell^2({\bb Z}_+) \subseteq \ell^2({\bb Z})$, where the latter is the Hilbert space of two-sided infinite square summable sequences, then using Fourier series we can identify the latter space with $L^2({\bb T})$, the space of square summable complex-valued measurable functions on the unit circle ${\bb T} = \partial{\bb D}$ with ${\bb D}$ the unit disk in ${\bb C}$. In this way, $\ell^2({\bb Z}_+)$ corresponds to the Hardy space $H^2({\bb D})$ of functions in $L^2({\bb T})$ which have holomorphic extensions to ${\bb D}$. The operator $U_+$ corresponds to the Toeplitz operator $T_z$ on $H^2({\bb D})$ defined to be multiplication by the coordinate function $z$. Moreover, the function $k_\lambda$ above corresponds to the analytic function $k_\lambda(z) = (1-\bar\lambda z)^{-1}$ for $z$ in ${\bb D}$.

If $\langle k_\lambda\rangle$ denotes the one-dimensional subspace of $H^2({\bb D})$ spanned by $k_\lambda$ for $\lambda$ in ${\bb D}$, then $\pi(\lambda) = \langle k_\lambda\rangle$ defines an anti-holomorphic function from ${\bb D}$ to the Grassmannian $Gr(1, H^2({\bb D}))$ of one-dimensional subspaces of $H^2({\bb D})$. This Grassmannian is an infinite dimensional complex manifold over which there is a canonical Hermitian holomorphic line bundle with the fiber over a point being the one-dimensional Hilbert subspace itself. We are interested in the natural pull-back bundle $E^*_{H^2({\bb D})}$ over ${\bb D}$ defined by $\pi$. (The presence of the asterisk corresponds to the fact that this bundle is anti-holomorphic. We see its dual, $E_{H^2({\bb D})}$, later in these notes.) Our goal is to study $H^2({\bb D})$ and the action of $T_z$ on it using $E^*_{H^2({\bb D})}$, an approach introduced by M.\ Cowen and the author in \cite{C-D}.

There are many other natural examples analogous to the Hardy space with similar structure. We describe one class, the family of weighted Bergman spaces on ${\bb D}$. Fix $\alpha$, $-1<\alpha<\infty$, and consider the measure $dA_\alpha = (1-|z|^2)^\alpha dA$ on ${\bb D}$, where $dA$ denotes area measure on ${\bb D}$. The weighted Bergman space $L^{2,\alpha}_a({\bb D})$ is the closed subspace of $L^2(dA_\alpha)$ generated by the polynomials ${\bb C}[z]$ or, equivalently, the functions in $L^2(dA_\alpha)$ which are equal  a.e.\ to a holomorphic function on ${\bb D}$.

There is an explicit formula for an analytic function $\gamma^\alpha_\lambda$ in $L^{2,\alpha}_a({\bb D})$ for $\lambda$ in ${\bb D}$ such that $M^*_{z,\alpha}\gamma^\alpha_\lambda = \bar\lambda \gamma^\alpha_\lambda$, where $M_{z,\alpha}$ is the operator on $L^{2,\alpha}_a({\bb D})$ defined to be multiplication by $z$ and such that $\gamma^\alpha_\lambda$ is anti-holomorphic in $\lambda$. In particular, one has
\begin{equation}\label{eq1.1}
\gamma^\alpha_\lambda(z) = (1-z\bar\lambda)^{-2-\alpha} = \sum^\infty_{n=0} \frac{\Gamma(n+2+\alpha)}{n!\Gamma(2+\alpha)}(z\bar\lambda)^n,
\end{equation}
where $\Gamma$ is the gamma function.
(For more information on weighted Bergman spaces, see \cite{HKZ}.) In particular, for $\alpha=0$, one obtains the classical Bergman space for ${\bb D}$, the closure of ${\bb C}[z]$ in $L^2(dA)$. 

Again one has an anti-holomorphic map $\pi_\alpha\colon \ {\bb D}\to Gr(1, L^{2,\alpha}_a({\bb D}))$ and a Hermitian anti-holomorphic line bundle $E^*_{L^{2,\alpha}_a({\bb D)})}$ over ${\bb D}$ for $-1<\alpha<\infty$.

\section{Quasi-Free Hilbert Modules}\label{sec1.2}

We have indicated that we are interested in studying not only the Hilbert space but also the action of certain natural operators on it. The best way to express this structure is in the language of Hilbert modules. Given an  $n$-tuple of commuting bounded operators $\pmb{T} = (T_1,\ldots, T_n)$ acting on a Hilbert space ${\cl H}$, there is a natural and obvious way of making ${\cl H}$ into a unital module over the algebra ${\bb C}[\pmb{z}] = {\bb C}[z_1,\ldots, z_n]$ of complex polynomials in $n$ commuting variables. In some cases, we are interested in contractive modules or module actions that satisfy $\|M_pf\|_{\cl H} \le \|p\|_{A({\bb B}^n)}\|f\|_{\cl H}$, where $M_p$ denotes the operator on ${\cl H}$ defined by module multiplication $M_pf = p\cdot f$. Here, $A({\bb B}^n)$ is the algebra of continuous functions on the closed unit ball in ${\bb C}^n$ that are holomorphic on ${\bb B}^n$, with the supremum norm.

More generally, we consider a unital module action of the function algebra $A(\Omega)$ on a Hilbert space ${\cl H}$, where $\Omega$ is a bounded domain in ${\bb C}^n$. Here $A(\Omega)$ is the closure in the supremum norm of the functions holomorphic on a neighborhood of the closure of $\Omega$.

A powerful technique from algebra, closely related to spectral theory, is localization. We sketch those aspects of this approach which we will need. For $\pmb{\omega}_0$ in ${\bb C}^n$, let
\[
I_{\pmb{\omega}_0} = \{p(\pmb{z}) \in {\bb C}[\pmb{z}]\colon \ p(\pmb{\omega}_0) = 0\}
\]
be the maximal ideal in ${\bb C}[\pmb{z}]$ of polynomials that vanish at $\pmb{\omega}_0$.

Now assume that ${\cl H}$ is a Hilbert module over ${\bb C}[\pmb{z}]$ for which
\begin{equation}\label{eq1.2}
\begin{split}
\text{(i)}~~&\text{DIM}_{\bb C} {\cl H}/I_{\pmb{\omega}}\cdot{\cl H} = m<\infty \text{ for  $\pmb{\omega}$ in ${\bb B}^n$ and}\\
\text{(ii)}~~&\bigcap\limits^\infty_{k=0} I^k_{\pmb{\omega}} \cdot {\cl H} = (0) \text{ for $\pmb{\omega}$ in ${\bb B}^n$.}
\end{split}
\end{equation}

Assumption (i) is usually referred to as the module  being semi-Fredholm (cf.\ \cite{E-P}). Note that the quotient module being finite dimensional implies that the submodule $I_{\pmb{\omega}}\cdot {\cl H}$ is closed in ${\cl H}$. The integer $m$ is called the multiplicity of ${\cl H}$. In general, this multiplicity is less than the rank of ${\cl H}$ or the smallest cardinality of a set of generators for ${\cl H}$ as a Hilbert module.

Hilbert modules satisfying (i) and (ii) have the properties referred to as quasi-freeness which we'll discuss later in these lectures. In \cite{DM2}, \cite{DM1}, properties of a closely related class of modules are obtained as well as alternate descriptions of the class. The concept is related to earlier work of Curto and Salinas \cite{C-S}. These modules  are viewed as the potential building blocks for general Hilbert modules. We will see in the second lecture that other sets of assumptions yield much of the structure possessed by  quasi-free Hilbert modules. In the following section we'll discuss  the intrinsic properties of the class. Finally, while we confine our attention here to Hilbert modules over the unit ball or over ${\bb C}[\pmb{z}]$, other domains in ${\bb C}^n$ are important.

There are many natural examples of such quasi-free Hilbert modules. In particular, if $L^2(\partial{\bb B}^n)$ denotes the usual Lebesgue space for surface measure on the unit sphere, $\partial {\bb B}^n$, then the Hardy module $H^2({\bb B}^n)$ can be defined as the closure of ${\bb C}[\pmb{z}]$ in $L^2(\partial {\bb B}^n)$. Similarly, the closure of ${\bb C}[\pmb{z}]$ in the Lebesgue space, $L^2({\bb B}^n)$, relative to volume measure on ${\bb B}^n$ yields the Bergman module $L^2_a({\bb B}^n)$. Both $H^2({\bb B}^n)$ and $L^2_a({\bb B}^n)$ are quasi-free Hilbert modules of multiplicity one. A related example is the $n$-shift module $H^2_n$, recently studied by several authors \cite{A1}, \cite{Dru}. One quick description of $H^2_n$ is that it is the symmetric Fock space in $n$ variables. While module multiplication by polynomials in ${\bb C}[\pmb{z}]$ define bounded operators on $H^2_n$, that is not the case for general functions in $A({\bb B}^n)$. However, $H^2_n$ is quasi-free having multiplicity one.

One could consider analogues of weighted Bergman spaces on ${\bb B}^n$, but we confine our attention here to these three examples.

The relation of quasi-free Hilbert modules to our earlier discussion of eigenvectors is given in the following lemma which is  straightforward to establish.

\begin{lem}\label{lem1.1}
A vector $f$ in the quasi-free Hilbert module ${\cl H}$ is orthogonal to $I_{\pmb{\omega}}\cdot {\cl H}$ for $\pmb{\omega}$ in ${\bb B}^n$ iff $M^*_pf = \ovl{p(\pmb{\omega})}f$ for $p$ in ${\bb C}[\pmb{z}]$.
\end{lem}

Thus ${\bb B}^n$ consists of common eigenvalues for the operators defined by the adjoint of module multiplication. The fact that one can find an anti-holomorphic function consisting of eigenvectors is taken up in the next section.

\section{Hermitian Holomorphic Vector Bundles}\label{sec1.3}

Let ${\cl H}_{\pmb{\omega}}$ denote the quotient Hilbert module ${\cl H}/I_{\pmb{\omega}} \cdot {\cl H}$ for $\pmb{\omega}$ in ${\bb B}^n$. Now the role of the complex conjugate is revealed since the map $\pmb{\omega}\to {\cl H}_{\pmb{\omega}}$ is ``holomorphic.'' However, this relationship is expressed not in terms of functions
but as sections of a holomorphic bundle.

For $f$ in ${\cl H}$, let $\hat f$ be the function in ${\bb B}^n$ so that $\hat f(\pmb{\omega})$ is the image of $f$ in ${\cl H}_{\pmb{\omega}}$.

\begin{defn}\label{defn1}
A Hilbert module ${\cl R}$ over ${\bb B}^n$ is said to be quasi-free of multiplicity $m$ if $\text{DIM}_{\bb C} {\cl R}_{\pmb{\omega}} = m$ for $\pmb{\omega}$ in ${\bb B}^n$ and $\bigcup\limits_{\pmb{\omega}\in {\bb B}^n} {\cl R}_{\pmb{\omega}}$ can be given the structure of a Hermitian holomorphic vector bundle over ${\bb B}^n$ so that the section $\hat f$ is holomorphic for $f$ in ${\cl R}$ and the map $f\to\hat f$ is one-to-one.
\end{defn}
 
We say that ${\cl R}$ is weakly quasi-free if these properties hold for ${\cl R}_{\pmb{\omega}} = {\cl R}/[I_{\pmb{\omega}}\cdot{\cl R}]$, where [~~] denotes closure.

\begin{thm}\label{thm1.1}
If ${\cl H}$ is a  Hilbert module satisfying conditions \eqref{eq1.2} for $m$ and $\pmb{\omega}$ in ${\bb B}^n$, then it is quasi-free having multiplicity $m$.

\end{thm}

This result is related to earlier work of Curto and Salinas \cite{C-S} and the existence of the bundle structure is established in Section \ref{sec2.2} \cite{DM2} with the additional assumption that ${\cl H}$ has a set of precisely $m$ generators as a Hilbert module. The  argument for the general case requires the use of Banach space-valued sheaf theory \cite{Put} and was described to the author by Putinar. However, the proof seems to require the assumption that ${\cl H}$ is finitely generated as a Hilbert module.\medskip 

\n {\bf Question 1.} Do conditions \eqref{eq1.2} imply that the theorem is true, at least for Hilbert modules over ${\bb B}^n$.

\begin{rem}\label{rem1.1}
While the sheaf-theoretic techniques used require that the module be finitely generated, it seems possible that one can restrict attention to open subsets $\Delta$ of ${\bb B}^n$ over which $\bigcup\limits_{\pmb{\omega}\in\Delta} {\cl R}_{\pmb{\omega}}$ is finitely generated. In particular, the same would be true on the overlap of two such open sets. This approach might allow one to prove a local version of Theorem \ref{thm1.1} and show that $\bigcup\limits_{\pmb{\omega}\in U} {\cl R}_{\pmb{\omega}}$ could be given the structure of a Hermitian holomorphic vector bundle, when the dimension of ${\cl R}_{\pmb{\omega}}$ is constant for $\pmb{\omega}$ in $U$ an open subset of ${\bb B}^n$. Moreover, it would follow that the map from ${\cl R}$ to holomorphic sections over $U$ would be injective.
\end{rem}

As the notation suggests, the bundle $E_{\cl H}$ is dual to the bundle $E^*_{\cl H}$ defined earlier (for the case $m=1$) as the pull-back of the anti-holomorphic map from ${\bb B}^n$ to $Gr(m,{\cl H})$ defined by $\pmb{\omega}\to (I_{\pmb{\omega}}\cdot {\cl H})^\bot$. We consider this result in the next lecture.

The importance of this representation theorem lies in the fact that the Hilbert module ${\cl H}$ is characterized by the Hermitian holomorphic vector bundle $E_{\cl H}$ and vice versa. This fact was established in the case $n=1$ or 2 by M.\ Cowen and the author in \cite{C-D}, \cite{C-D2}; for commuting $n$ tuples, $n>1$, by Curto and Salinas \cite{C-S}; and for Hilbert modules by Chen and the author \cite{Ch-D}.

\begin{thm}\label{thm1.2}
Two quasi-free Hilbert modules ${\cl H}$ and $\widetilde{\cl H}$ over ${\bb B}^n$, are isomorphic iff the corresponding Hermitian holomorphic vector bundles $E_{\cl H}$ and $E_{\widetilde{\cl H}}$ are isomorphic.
\end{thm}

More precisely, there exists a unitary module map $U\colon \ {\cl H}\to \widetilde{\cl H}$ iff there exists an isometric holomorphic bundle map $\Phi\colon \ E_{\cl H}\to E_{\widetilde{\cl H}}$.

That $U$ is a module map means that $UM_p = \widetilde M_pU$ for $p$ in ${\bb C}[\pmb{z}]$, where $\widetilde M_p$ denotes the operators defined by module multiplication on $\widetilde{\cl H}$. Thus $U$ identifies ${\cl H}$ and $\widetilde{\cl H}$ as Hilbert modules.

The statement about $\Phi$ means that the bundle map is holomorphic from the total space of $E_{\cl H}$ to that of $E_{\widetilde{\cl H}}$ and acts as a unitary operator from the fiber $E_{\cl H}|_{\pmb{\omega}}$ to the corresponding fiber  $E_{\widetilde{\cl H}}|_{\pmb{\omega}}$ for $\pmb{\omega}$ in ${\bb B}^n$.

A key step in the proof of this result depends on a Rigidity Theorem, established in full generality in \cite{C-D}.

\begin{thm}\label{thm1.3}
If $\Omega$ is a domain in ${\bb C}^n$, $k$ a positive integer, ${\cl H}$ and $\widetilde{\cl H}$ are Hilbert spaces with $\pi$ and $\tilde\pi$ anti-holomorphic maps from $\Omega$ to $Gr(k,{\cl H})$ and $Gr(k,\widetilde{\cl H})$, respectively, then the Hermitian anti-holomorphic pull-back bundles for $\pi^*$ and $\tilde\pi^*$ over $\Omega$ are equivalent iff there exists a unitary $U\colon \ {\cl H}\to \ovl{\cl H}$ such that
\[
\pi(\pmb{\omega}) = U^*\tilde\pi(\pmb{\omega}) U \quad \text{for $\omega$ in } \Omega.
\]
\end{thm}

An analogous result holds for the quotient Grassmanian which allows one to replace $E^*_{\cl H}$ and $E^*_{\widetilde{\cl H}}$ in Theorem \ref{thm1.2} by $E_{\cl H}$ and $E_{\widetilde {\cl H}}$.

\section{Curvature- and Operator-Theoretic Invariants}\label{sec1.4}

Now Theorem \ref{thm1.2} is particularly useful in this context because we can directly relate operator-theoretic invariants with those from complex geometry. Recall that in the middle of the twentieth century, Chern observed that there is a unique canonical connection on a Hermitian holomorphic bundle (cf.\ \cite{Wells}). A connection is a first order differential operator defined on the smooth sections of the vector bundle. Connections always exist and can be expressed in terms of a sum of smooth one forms. In the holomorphic context, the connection can be expressed as the sum of (1,0) forms and (0,1) forms. The key observation of Chern is that when the bundle is Hermitian and holomorphic  one can require that the (0,1) forms all be zero and in this case the connection is  unique.  The curvature ${\cl K}(\pmb{\omega})$ for a connection is a section of the two-forms built on the bundle. Hence, one can speak of \emph{the} curvature of a Hermitian holomorphic vector bundle and hence, in particular, in our context. In \cite{C-D}, Cowen and the author showed how to calculate this curvature and its partial derivatives in terms of  operator theoretic invariants. 

Let us assume that ${\cl H}$ is a quasi-free, multiplicity one Hilbert module over ${\bb D}$. Then for $\omega$ in ${\bb D}$, $\text{NULL}(M_z-\omega)^{*2}$ is two-dimensional and the restriction of $(M_z-\omega)^*$ to it is a nilpotent operator of order two. If one chooses an orthonormal basis for $\text{NULL}(M_z-\omega)^{*2}$ correctly, then $(M_z-\omega)^*$ has the form
\[
\begin{pmatrix}
0&h_{\cl H}(\omega)\\ 0&0
\end{pmatrix}.
\]
If one requires that $h_{\cl H}(\omega)>0$, then the function $h_{\cl H}(\omega)$ is unique. Further, one can show that ${\cl K}_{\cl H}(\omega) = -\frac1{h_{\cl H}(\omega)^2} d\omega d\bar\omega$ in general. Since the curvature is a complete invariant for line bundles, we see that this function $h_{\cl H}(\omega)$ is a complete invariant for the Hilbert module ${\cl H}$.

We can calculate this function for the examples introduced earlier. In particular,
\begin{equation}\label{eq1.3}
\begin{split}
h_{H^2({\bb D})}(\omega) &= 1 -|\omega|^2, h_{L^2_a({\bb D})}(\omega) = \frac1{\sqrt 2} (1-|\omega|^2)\quad \text{and}\\
h_{L^{2,\alpha}_a({\bb D})}(\omega) &= \frac1{\sqrt{2+\alpha}} (1-|\omega|^2).
\end{split}
\end{equation}
These calculations yield the fact that none of these Hilbert modules are isomorphic. (Of course, there are other ways to see this fact.)

One can also calculate the curvatures for the Hardy and Bergman modules for ${\bb B}^n$ but curvature in this case can not be described by a single function. We'll say more about this case in the next lecture.

\section{Reducing Submodules}\label{sec1.5}

We continue by relating some operator theoretic concepts to their complex geometric counter-parts. Since the latter invariants form a complete set, in principle, we can always do that but in some cases the results are particularly striking and useful.

We begin by considering reducing submodules. For general quasi-free Hilbert modules, there is no characterization available for general submodules. Consideration of the situation for just the Hardy and Bergman modules over ${\bb D}$  shows just how complex is the structure. Still we will have some things to say about submodules  in Section \ref{sec2.5}. Here we want to consider submodules ${\cl L}_1$ contained in a quasi-free Hilbert module ${\cl H}$ for which there exists a submodule ${\cl L}_2$ such that ${\cl H} = {\cl L}_1\oplus {\cl L}_2$. (The symbol $\oplus$ indicates an orthogonal direct sum.) These are the \emph{reducing submodules}.

A standard argument shows that a reducing submodule is the range of a projection $P$ which is a module map. Hence, one approach to characterizing reducing submodules is to consider first the commutant of the Hilbert module; that is, all $X$ in ${\cl L}({\cl H})$ which are module maps. We begin with the case of a quasi-free Hilbert module ${\cl R}$ over ${\bb B}^n$ having multiplicity one. If $X\colon \ {\cl R}\to {\cl R}$ is a module map, then we can define the complex number $\widehat X(\pmb{\omega})$ so that $\widehat X(\pmb{\omega}) \hat f(\pmb{\omega}) = (\widehat{Xf})(\pmb{\omega})$ for $\pmb{\omega}$ in ${\bb B}^n$ and $f$ in ${\cl H}$. It is easy to see that $|\widehat X(\omega)| \le \|X\|$ and hence $\widehat X$ is in the algebra, $H^\infty({\bb B}^n)$, of bounded holomorphic functions on ${\bb B}^n$. For $H^2({\bb B}^n)$, $L^{2,\alpha}_a({\bb D})$ and $L^2_a({\bb B}^n)$, one can show that the commutant equals $H^\infty({\bb B}^n)$. This is not the case in general. In particular, it is not true for $H^2_n$ for $n>1$ \cite{A1}.

If we consider a quasi-free Hilbert module ${\cl R}$ having multiplicity $m>1$, then $\widehat X$ is a bounded holomorphic bundle map on $E_{\cl R}$ and should not be thought of as a function.

We  now return to the question of characterizing reducing submodules.
If $X$ is a projection in the commutant of ${\cl R}$, then $\widehat X(\pmb{\omega})$ is a projection on the fiber $E_{\cl R}|_{\pmb{\omega}}$. If we set $F_{\pmb{\omega}} = \text{RAN}\{\widehat X(\pmb{\omega})\}$ and $\widetilde F_{\pmb{\omega}} = \text{RAN}\{(\widehat{I-X})(\pmb{\omega})\}$ for $\pmb{\omega}$ in ${\bb B}^n$, then we obtain orthogonal holomorphic sub-bundles $F$ and $\widetilde F$ of $E_{\cl R}$. The converse is true which we will prove in Section \ref{sec2.3}.

\begin{thm}\label{thm1.4}
If ${\cl R}$ is a quasi-free Hilbert module of finite multiplicity $m$ over ${\bb B}^n$, then ${\cl L}$ is a reducing submodule of ${\cl R}$ iff there exists  holomorphic sub-bundles $F$ and $\widetilde F$ of $E_{\cl R}$ so that $E_{\cl R} = F\oplus \widetilde F$ and ${\cl L} = \{f\in{\cl R}\mid \hat f(\pmb{\omega})\in F_{\pmb{\omega}} ~\forall \pmb{\omega} \in {\bb B}^n\}$.
\end{thm}

\begin{cor}
If ${\cl R}$ is a quasi-free Hilbert module over ${\bb B}^n$ of multiplicity one, then ${\cl R}$ is irreducible.
\end{cor}

These arguments extend to the context of complementary submodules since such a submodule ${\cl M}$ of the quasi-free Hilbert module ${\cl R}$ is the range of an idempotent map in the commutant of ${\cl R}$. Thus the same argument establishes one direction of the analogue of the previous theorem to obtain:

\begin{thm}\label{thm1.5}
If ${\cl R}$ is a quasi-free Hilbert module of finite multiplicity over ${\bb B}^N$ and ${\cl L}$ is a complemented submodule, then there exists a holomorphic sub-bundle $F$ of $E_{\cl R}$ such that ${\cl L} = \{f\in {\cl R}\mid \hat f(\omega) \in F_{\pmb{\omega}} \forall \pmb{\omega} \in {\bb B}^n\}$.
\end{thm}

The converse is not clear. Since there is a necessary condition involving the norm of the localization of the idempotent map. Whether this assumption is sufficient is unknown.
\medskip 

\n {\bf Question 2.} Let ${\cl R}$ be a quasi-free Hilbert module of finite multiplicity over ${\bb B}^n$ so that there exist holomorphic sub-bundles $F_1$ and $F_2$ with $E_{\cl R}|_{\pmb{\omega}} = F_1|_{\pmb{\omega}} \dot+ F_2|_{\pmb{\omega}}$ for $\pmb{\omega}$ in ${\bb B}^n$. If the angle between $F_1|_{\pmb{\omega}}$ and $F_2|_{\pmb{\omega}}$ is bounded away from 0, does there exist submodules ${\cl L}_1$ and ${\cl L}_2$ related to $F_1$ and $F_2$ as in Theorem \ref{thm1.5} so that ${\cl R} = {\cl L}_1 \dot+ {\cl L}_2$?

The preceding results raise other interesting questions about the situation for submodules.\medskip

\n {\bf Question 3.} Let ${\cl R}$ be a quasi-free Hilbert module over ${\bb B}^n$ of multiplicity one and ${\cl L}$ be a submodule of ${\cl R}$. Does it follow that ${\cl L}$ is irreducible? If so, does ${\cl L}$ have any complemented submodules? What if ${\cl L}$ is finitely generated?

\begin{rem}\label{rem1.2}
The answer is unknown for even weaker questions such as the following:\ Do there exist non-zero submodules ${\cl L}_1$ and ${\cl L}_2$ of ${\cl R}$ is so that ${\cl L}_1 \perp {\cl L}_2$? This is not possible if ${\cl R}$ is subnormal which is defined in Section \ref{sec3.1}.
\end{rem}

\begin{rem}\label{rem1.3}
Since Rudin exhibited \cite{Rud} submodules of the Hardy module over the bidisk which are not finitely generated, it seems likely that the same might be true for quasi-free Hilbert modules over ${\bb B}^n$.
\end{rem}

The relevance of this finiteness assumption arises when one considers the ``spectral sheaf'' $\bigcup\limits_{\pmb{\omega} \in {\bb B}^n} {\cl L}_{\pmb{\omega}}$ over ${\bb B}^n$. If there exists $\pmb{\omega}_0$ in ${\bb B}^n$ for which $\text{DIM}_{\bb C} {\cl L}_{\pmb{\omega}_0} = 1$, then following up on the suggestion in Remark \ref{rem1.1}, it might be possible to make $\bigcup\limits_{\pmb{\omega}\in\Delta} {\cl L}_{\pmb{\omega}}$ into a Hermitian holomorphic line bundle and use the same proof as that used for Theorem \ref{thm1.1}. It seems probable that for finitely generated ${\cl L}$, $\text{DIM}_{\bb C} {\cl L}_{\pmb{\omega}} \ge 2$ for all $\pmb{\omega}$ in ${\bb B}^n$ is not possible, but this is unknown.

\section{Quasi-Similar Hilbert Modules}\label{sec1.6}

Analogous to the commutant for one Hilbert module is the space of intertwining module maps between two quasi-free Hilbert modules. Again such a map leads to a map between the corresponding bundles. Utilizing this object, however, can be somewhat difficult. For example, the question of deciding when two modules are similar or quasi-similar comes down to determine the existence or non-existence of such maps, a problem which hasn't been solved in general. But sometimes the difficulties can be overcome as was observed in the following case by Curto and Salinas \cite{C-S}.

Suppose $X\colon \ L^{2,\alpha}_a({\bb D}) \to L^{2,\beta}_a({\bb D})$ is an intertwining module map for $-1<\alpha,\beta < \infty$ and let $\gamma^\alpha_\omega$ and $\gamma^\beta_\omega$ be the anti-holomorphic functions defined in Section \ref{sec1.1} which satisfy $M^*_z\gamma^\alpha_\omega = \bar\omega\gamma^\alpha_\omega$ and $M^*_z\gamma^\beta_\omega = \bar\omega \gamma^\beta_\omega$ for $\omega$ in ${\bb D}$.
Then there exists an anti-holomorphic function $\varphi$ on ${\bb D}$ such that $X^*\gamma^\beta_\omega = \varphi(\omega)\gamma^\alpha_\omega$ and hence
\begin{equation}\label{eq1.4}
|\varphi(\omega)| \le \|X^*\| \frac{\|\gamma^\beta_\omega\|}{\|\gamma^\alpha_\omega\|} = \|X^*\| (1-|\omega|^2)^{\alpha-\beta}.
\end{equation}
Thus, if $\alpha-\beta>0$, we have $\varphi\equiv 0$ and hence $X=0$.

\begin{thm}\label{thm1.6}
For $-1<\alpha\ne \beta<\infty$, the Hilbert modules $L^{2,\alpha}_a({\bb D})$ and $L^{2,\beta}_a({\bb D})$ are not quasi-similar; that is, there do not exist module maps $X,Y$, $X\colon \ L^{2,\alpha}_a({\bb D})\to L^{2,\beta}_a({\bb D})$ and $Y\colon \ L^{2,\beta}_a({\bb D}) \to L^{2,\alpha}_a({\bb D})$
 which are injective with dense range.  
\end{thm}

Actually the preceding proof shows that there is no non-zero module map $X\colon \ L^{2,\alpha}_a({\bb D}) \to L^{2,\beta}_a({\bb D})$ when $\alpha<\beta$.

Although much effort has been directed to understanding when two quasi-free Hilbert modules of finite multiplicity are similar (cf.\ \cite{JW}), very few results have been obtained, at least based on complex geometry.
\medskip 

\n {\bf Question 4.} Can one give conditions involving the curvatures which imply that two quasi-free Hilbert modules of multiplicity one are similar?

\clearpage\thispagestyle{empty}\cleardoublepage

\chapter{SECOND LECTURE}\label{chp2}

\section{Module Resolutions}\label{sec2.1}

A fundamental approach to the study of contraction operators on Hilbert space is the model theory of Sz-Nagy and Foias \cite{Sz-NF}. Let us recall some aspects of it, placed in the context of contractive Hilbert modules over the disk algebra $A({\bb D})$.

First, recall the von Neumann inequality states that $\|p(T)\|_{{\cl L}({\cl H}_T)} \le \|p\|_{A({\bb D})}$ for $T$ a contraction on the Hilbert space ${\cl H}_T$ and $p(z)$ in ${\bb C}[z]$. Using this result one sees that the contraction $T$ on ${\cl H}_T$ can be used to make ${\cl H}_T$ into a contractive Hilbert module over $A({\bb D})$.
And, vice versa, module multiplication by the coordinate function $z$ recovers the contraction operator $T$ on ${\cl H}_T$.

Second, if for the contraction operator $T$, the sequence $T^n$ converges strongly to 0 (or $T$ belongs to class $C_{0\cdot}$), then there exist coefficient Hilbert spaces ${\cl D}$ and ${\cl D}_*$ and module maps $X\colon \ H^2_{\cl D}({\bb D})\to {\cl H}_T$ and $Y\colon\ H^2_{{\cl D}_*}({\bb D})\to H^2_{\cl D}({\bb D})$ such that $X$ is a co-isometry, $Y$ is an isometry, and the sequence
\begin{equation}\label{eq2.1}
0 \longleftarrow {\cl H}_T \overset{\sst X}{\longleftarrow} H^2_{\cl D}({\bb D}) \overset{\sst Y}{\longleftarrow} H^2_{{\cl D}_*}({\bb D}) \longleftarrow 0 \quad \text{ is exact.}
\end{equation}
Since $X$ is onto and $Y$ is one-to-one, exactness requires only the additional condition  that the range of $Y$ equals the null space of $X$. If one localizes the map $Y$ over ${\bb D}$; that is, one considers the map
\begin{equation}\label{eq2.2}
H^2_{\cl D}({\bb D})/I_\omega\cdot H^2_{\cl D}({\bb D}) \overset{\sst Y_\omega}{\longleftarrow} H^2_{{\cl D}_*}({\bb D})/I_\omega \cdot H^2_{{\cl D}_*}({\bb D})
\end{equation}
defined for $\omega$ in ${\bb D}$, one obtains the characteristic operator function $\Theta_T(\omega)$ for the operator $T$. An important ingredient in the theory is an explicit formula for $\Theta_T(\omega)$ in terms of $T$ and ${\cl H}_T$. In particular, one takes ${\cl D}$ and ${\cl D}_*$ to be  the closures of the ranges of $D_T = (I-T^*T)^{\frac12}$ and $D_{T^*} = (I-TT^*)^{\frac12}$, respectively, and $\Theta_T(z)\colon \ {\cl D}_T \to {\cl D}_{T^*}$, where
\begin{equation}\label{eq2.3}
\Theta_T(z) = [-T+zD_{T^*}(I-zT^*)^{-1}]|_{{\cl D}_T}.
\end{equation}

A particularly simple case of this formula occurs for Jordan models of multiplicity one which are defined by inner functions. (Recall that an \emph{inner function} is a function $\theta$ in $H^\infty({\bb D})$ with unimodular non-tangential boundary values a.e.\ on $\partial {\bb D}$.) More specifically, if ${\cl H}_\theta$ is the quotient module $H^2({\bb D})/\theta H^2({\bb D})$ defined by the inner function $\theta$, then we have the resolution
\begin{equation}\label{eq2.4}
0 \longleftarrow {\cl H}_{\theta} \overset{\sst X}{\longleftarrow} H^2({\bb D}) \overset{\sst Y}{\longleftarrow} H^2({\bb D}) \longleftarrow 0
\end{equation}
with $Y$ being defined to be multiplication by $\theta$. Localization of the module map $Y$ yields the fact that $\widehat Y(z) = \theta(z)$ which is equivalent to \eqref{eq2.3}. (There is an arbitrary scalar of modulus one here which arises from the non-canonical identification of the line bundles determined by the two Hardy modules.)

In the introductory presentation on Hilbert modules given by Paulsen and the author \cite{D-P}, this interpretation of the canonical model was given noting the analogue of this resolution of ${\cl H}_T$ in terms of $H^2_{\cl D}({\bb D})$ and $H^2_{{\cl D}_*}({\bb D})$ with the projective resolutions used in algebra to study more general modules. The goal posed in \cite{D-P} was to construct resolutions for a Hilbert module in terms of  ``\v Silov modules'' (which we will discuss in the third lecture) for contractive Hilbert modules over function algebras. Since those notes were written, it has become clear that this is not the best approach. Rather the author now believes the key to the effectiveness of the Sz.-Nagy--Foias model theory  rests, in large part, on the fact that (1) the Hardy modules, $H^2_{\cl D}({\bb D})$ and $H^2_{{\cl D}_*}({\bb D})$, are quasi-free in that $H^2_{\cl D}({\bb D})/I_\omega\cdot H^2_{\cl D}({\bb D})$ and $H^2_{{\cl D}_*}({\bb D})/I_\omega\cdot H^2_{{\cl D}_*}({\bb D})$ are isomorphic to ${\cl D}$ and ${\cl D}_*$, respectively, for $\omega$ in ${\bb D}$ and (2) these modules give rise to Hermitian holomorphic vector bundles over ${\bb D}$ with the characteristic operator function  being a holomorphic bundle map. In any case,  the goal we are currently pursuing is to construct resolutions for Hilbert modules in terms of quasi-free Hilbert module (cf.\ \cite{DM2}, \cite{DM1}). In these notes, we will say little more  about such resolutions except for the fact that this goal provide much of our motivation for studying quasi-free Hilbert modules. Of course, there is also the fact that the class of quasi-free Hilbert modules contains most of the classical examples of Hilbert spaces of holomorphic functions closed under multiplication by polynomials.

There is another technique important in algebraic geometry which we are mimicking in our approach, that of the resolution of sheaves by vector bundles. More specifically, localization of a general Hilbert module ${\cl M}$ yields a notion of spectral sheaf $\bigcup\limits_{\pmb{\omega}\in {\bb  B}^n} {\cl M}_{\pmb{\omega}}$ in many cases. (This spectral sheaf has few nice properties unless we make assumptions about ${\cl M}$.) The resolution by quasi-free Hilbert modules described above should provide a resolution of this sheaf in terms of holomorphic vector bundles. While there is a wealth of technicalities to overcome to validate this approach, it is one picture we have in mind for the future development of multivariate operator theory.

\section{Quasi-Free Hilbert Modules Revisited}\label{sec2.2}

While we introduced quasi-free Hilbert modules in the first lecture, we return to the topic with  an alternate approach as well as enough details to make clear the meaning of the statement in the first lecture that the bundles $E_{\cl H}$ and $E^*_{\cl H}$ are dual to each other.

Although the spaces ${\cl D}$ and ${\cl D}_*$ are, in general, infinite dimensional, the model theory is particularly powerful when they are finite dimensional and $\Theta_T$ is essentially a  holomorphic matrix-valued function. As in the first lecture, we will confine our attention to quasi-free Hilbert modules of finite multiplicity. Although such modules are not sufficient to form the building blocks for general modules, the class that does possess such resolutions is likely to be an interesting one, whose study should reveal many interesting insights and results, particularly for multivariate operator theory. 

We continue to restrict attention to modules over the unit ball and the ball algebra $A({\bb B}^n)$. Free modules over $A({\bb B}^n)$ in the sense of algebra have the form $A({\bb B}^n) \otimes_{\text{alg}} {\bb C}^k$, or at least those which are finitely generated. (The subscript ``alg'' on $\otimes$ indicates that the symbol denotes the algebraic tensor product.) However, $A({\bb B}^n) \otimes_{\text{alg}} {\bb C}^k$ is \emph{NOT} a Hilbert space. Thus we consider inner products on $A({\bb B}^n) \otimes_{\text{alg}} {\bb C}^k$ and their completions. Not every inner product can be used since we want to preserve the holomorphic character of the elements of the completion. Thus we assume the inner product $\langle ~~,~~\rangle_{\cl R}$ satisfies:
\begin{align}\label{eq2.5}
\text{(1)}~~&\text{eval$_{\pmb{\omega}}\colon \ A({\bb B}^n) \otimes_{\text{alg}} {\bb C}^k \to {\bb C}^k$ is bounded and locally uniformly bounded}\notag\\
&\text{in the induced ${\cl R}$-norm;}\notag\\
\text{(2)}~~&\left\|\sum\limits^k_{i=1} p\varphi_i\otimes e_i\right\|_{\cl R} \le \|p\|_{A({\bb B}^n)} \left\|\sum\limits^\ell_{i=1} \varphi_i\otimes e_i\right\|_{\cl R} \text{ for $p$ in ${\bb C}[\pmb{z}]$, $\{\varphi_i\}$ in $A({\bb B}^n)$,}\notag\\
&\text{and $\{e_i\}$ in ${\bb C}^k$; and}\\
\intertext{\newpage}
\text{(3)}~~&\text{for $\left\{\sum\limits^k_{i=1} \varphi^{(\ell)}_i\otimes e_i\right\}$ Cauchy in the ${\cl R}$-norm, $\lim\limits_{\ell\to\infty} \sum\limits^k_{i=1} \varphi^{(\ell)}_i (\pmb{\omega}) \otimes e_i = 0$}\notag\\
&\text{for all $\pmb{\omega}$ in ${\bb B}^n$ iff}\notag
\end{align}
\[
\lim_{\ell\to\infty} \left\|\sum \varphi^{(\ell)}_i \otimes e_i\right\|_{\cl R} = 0.
\]
We let ${\cl R}$ denote the completion of $A({\bb B}^n) \otimes_{\text{alg}} {\bb C}^k$ in the ${\cl R}$-norm.

Condition (1) implies that the completion in the ${\cl R}$-norm can be identified with ${\bb C}^k$-valued functions on ${\bb B}^n$ and the local uniform boundedness implies that these functions are holomorphic. Condition (3) implies that the limit function of a Cauchy sequence in $A({\bb B}^n) \otimes_{\text{alg}} {\bb C}^k$ vanishes identically iff the limit in the ${\cl R}$-norm is the zero function. Finally, condition (2) ensures that ${\cl R}$ is a contractive Hilbert module on $A({\bb B}^n)$.

A Hilbert module ${\cl R}$ so obtained does not have the same properties as those described in \eqref{eq1.2}. In that case, ${\cl R}$ has $k$ generators as a Hilbert module over $A({\bb B}^n)$ and $\text{DIM}_{\bb C}{\cl R}/[I_{\pmb{\omega}}\cdot {\cl R}]_{\cl R} = k$ for $\pmb{\omega}$ in ${\bb B}^n$, where the bracket $[I_{\pmb{\omega}}\cdot {\cl R}]$ denotes closure in the ${\cl R}$ norm. Conditions \eqref{eq1.2} don't guarantee that ${\cl R}$ is finitely generated let alone $k$-generated. Moreover, in general, we can't conclude from \eqref{eq2.5} that $I_{\pmb{\omega}}\cdot {\cl R}$ is closed in ${\cl R}$. In \cite{DM2}, however, the bundle structure for $\bigcup\limits_{\pmb{\omega}\in {\bb B}^n} {\cl R}/[I_{\pmb{\omega}}\cdot{\cl R}]$ is demonstrated under conditions \eqref{eq2.5}, which we'll restate as a theorem.

\begin{thm}\label{thm2.1}
If ${\cl H}$ is a Hilbert module over $A({\bb B}^n)$ obtained as the completion of $A({\bb B}^n) \otimes_{\text{alg}} {\bb C}^k$ in a norm satisfying conditions \eqref{eq2.5}, then it is weakly quasi-free.
\end{thm}

We now want to consider the relationship between $E_{\cl H}$ and $E^*_{\cl H}$, for weakly quasi-free Hilbert modules under the assumption that conditions \eqref{eq2.5} hold.

If ${\cl R}$ is the completion of $A({\bb B}^n) \otimes_{\text{alg}} {\bb C}^m$ and $\{e_i\}^m_{i=1}$ is a basis for ${\bb C}^m$, then set $k_i=1 \otimes e_i$ for $i=1,2,\ldots, m$. Then $\{k_i\}$ is a set of module generators for ${\cl R}$. Using a standard argument, we can identify $f$ in ${\cl R}$ with a holomorphic function $\hat f\colon \ {\bb B}^n \to {\bb C}^m$ so that $\hat k_i(\pmb{\omega}) \equiv e_i$ for $\pmb{\omega}$ in ${\bb B}^n$ and $i=1,2,\ldots, m$. Note that we can identify ${\cl R}_{\pmb{\omega}} = {\cl R}/[I_{\pmb{\omega}}\cdot {\cl R}]$ canonically with ${\bb C}^m$ and hence we can identify $E_{\cl R} = \bigcup\limits_{\pmb{\omega}} {\cl R}_{\pmb{\omega}}$ with the trivial bundle ${\bb B}^n\times {\bb C}^m$ over ${\bb C}^m$. Moreover, the correspondence $f\to \hat f$ yields the desired injection of ${\cl R}$ into ${\cl O}(E_{\cl R}) \cong {\cl O}({\bb B}^n, {\bb C}^m)$, the spaces of holomorphic sections of $E_{\cl R}$ and holomorphic ${\bb C}^m$-valued functions on ${\bb B}^n$, respectively.

For each $\pmb{\omega}$ in ${\bb B}^n$ and $i=1,\ldots, m$, there exists a vector $h_i(\pmb{\omega})$ in ${\cl R}$ so that $\langle f,h_i(\pmb{\omega})\rangle_{\cl R} = \langle \hat f(\pmb{\omega}), e_i\rangle_{{\bb C}^m}$. The functions $\{h_i\}$ are anti-holomorphic and the $\{h_i(\pmb{\omega})\}$ forms a basis for $(I_{\pmb{\omega}} \cdot {\cl R})^\bot$. Hence, the $\{h_i\}$ forms an anti-holomorphic frame for the dual bundle $E^*_{\cl R}$. The duality between the bases $\{k_i(\pmb{\omega})\}^m_{i=1}$ and $\{h_i(\pmb{\omega})\}^m_{i=1}$ establishes the duality between the Hermitian holomorphic vector bundle $E_{\cl R}$ and the Hermitian anti-holomorphic vector bundle $E^*_{\cl R}$. A similar result holds for the quasi-free case, in general, but here one must work locally since one doesn't have a global frame for $E_{\cl R}$ consisting of sections defined by vectors in ${\cl R}$. One knows that the bundle $E_R$ does have a holomorphic frame since all holomorphic vector bundles over ${\bb B}^n$ are trivial. However, the key is whether such a frame can be obtained from sections defined by vectors in the Hilbert space.
\medskip

\n {\bf Question 5.} Does there always exist a frame for $E_{\cl R}$ determined by a finite number of elements of ${\cl R}$?

This question is equivalent to whether (weakly) quasi-free Hilbert modules of finite multiplicity are always finitely generated.

\section{Reducing Submodules Revisited}\label{sec2.3}

In the previous lecture, reducing submodules of a quasi-free Hilbert module were considered. In  the multiplicity one case, it was shown that such a module is irreducible and there are no reducing submodules. Thus the Hardy module and the weighted Bergman modules on the disk, or on the unit ball, are all irreducible. The same is true for $H^2_n$. Still, there is a family of quasi-free Hilbert modules which can be constructed from them for which the answer is less obvious and more interesting. However, before we consider this class, let us take a closer look at the reducing submodules of a quasi-free Hilbert module ${\cl R}$ and their relation to the Hermitian holomorphic vector bundle $E_{\cl R}$. We begin with a lemma.

\begin{lem}\label{lem2.1}
Let ${\cl R}$ be a weakly quasi-free Hilbert module over ${\bb B}^n$ of multiplicity $m$, $m<\infty$, $\pmb{\omega}_0$ in ${\bb B}^n$, $V$ a neighborhood of $\pmb{\omega}_0$ in ${\bb B}^n$, and $\{k_i\}^m_{i=1}$ vectors in ${\cl R}$ such that $\{\hat k_i(\pmb{\omega})\}^m_{i=1}$ forms a basis for ${\cl R}_{\pmb{\omega}}$ for $\pmb{\omega}$ in $V$. Then $\bigvee\limits_{\pmb{\alpha},i} \frac{\partial^{\pmb{\alpha}}}{\partial\bar{\pmb{z}}^{\pmb{\alpha}}} k_i(\pmb{\omega}_0) = {\cl R}$.
\end{lem}

Here $\pmb{\alpha}$ is the multi-index $(\alpha_1,\ldots, \alpha_n)$ with each $\alpha_i$ a non-negative integer.

The argument uses a vector $f$ in ${\cl R}$ to reduce the proof to the fact that the only vector-valued holomorphic function for which all partial derivatives vanish at a point is the zero function.

\begin{proof}[Proof of Theorem \ref{thm1.4}]
Suppose ${\cl L}$ is a reducing subspace of ${\cl R}$ and $P$ is the orthogonal projection of ${\cl R}$ onto ${\cl L}$. A straightforward argument shows that $P(I_{\pmb{\omega}}\cdot {\cl R})^\bot \subset (I_\omega\cdot {\cl R})^\bot$ and similarly, $P^\bot(I_{\pmb{\omega}}\cdot {\cl R})^\bot \subset (I_\omega\cdot{\cl R})^\bot$, where $P^\bot = I-P$. Thus each subspace $(I_{\pmb{\omega}}\cdot {\cl R})^\bot = P(I_{\pmb{\omega}}\cdot {\cl R})^\bot \oplus P^\bot(I_{\pmb{\omega}}\cdot {\cl R})^\bot$. This is, of course the decomposition of $E^*_{\cl R}$ into an orthogonal direct sum or as the direct sum of two Hermitian anti-holomorphic bundles. In particular, using the dual basis $\{h_i\}^m_{i=1}$, we obtain the sets $\{Ph_i\}^m_{i=1}$ and $\{P^\bot h_i\}^m_{i=1}$, which are spanning holomorphic sections for the two bundles. An easy argument involving dimension shows that these sections span the fiber but don't necessarily form a basis.

For the other direction of the argument, suppose we can write $E^*_{\cl R} = F_1 \oplus F_2$, where $F_1$ and $F_2$ are anti-holomorphic sub-bundles. (Here, the symbol $\oplus$ indicates that $F_1|_{\pmb{\omega}}\perp F_2|_{\pmb{\omega}}$ for $\pmb{\omega}$ in ${\bb B}^n$.) Take local anti-holomorphic sections $\{g_i\}^m_{i\equiv 1}$ of $E^*_{\cl R}$ such that $g_1,\ldots, g_{m_0}$ span $F_1$ and $\{g_{m_0+1},\ldots, g_m\}$ span $F_2$ in a neighborhood of a point $\pmb{\omega}_0$ in ${\bb B}^n$. Since $\langle g_i(\pmb{\omega}), g_j(\pmb{\omega})\rangle_{\cl R} = 0$ for $1\le i\le m_0$ and $m_0< j\le m$, we can differentiate by $\frac{\partial^{\pmb{\alpha}+\pmb{\beta}}}{\partial \pmb{z}^{\pmb{\alpha}} \partial \bar{\pmb{z}}^{\pmb{\beta}}}$ to show that $\bigvee\limits_{\underset{\sst 1\le i\le m_0}{\pmb{\alpha}}} \frac{\partial^{\pmb{\alpha}}}{\partial \bar{\pmb{z}}^{\alpha}} g_i(\pmb{\omega}_0)$ and $\bigvee\limits_{\underset{\sst m_0<j\le m}{\pmb{\beta}}} \frac{\partial^{\pmb{\beta}}}{\partial\bar{\pmb{z}}^{\pmb{\beta}}} g_j(\pmb{\omega}_0)$ are orthogonal subspaces of ${\cl R}$ which span ${\cl R}$ by the lemma. Therefore, the decomposition $E_{\cl R} = F_1\oplus F_2$ yields reducing submodules.
\end{proof}

One conclusion we can draw from this line of argument is that reducing submodules determine orthogonal decompositions of a fiber $E^*_{\cl H}|_\omega$ relative to which the curvature matrix also decomposes as an orthogonal direct sum. 

\begin{thm}\label{thm2.2}
Let ${\cl R}$ be a quasi-free Hilbert module over $A({\bb D})$ of finite multiplicity and $\omega_0$ be a point in ${\bb D}$ at which the curvature ${\cl K}_{\cl R} d\omega  d\bar\omega$ has distinct eigenvalues. Then the lattice of reducing submodules of ${\cl R}$ is finite and discrete.
\end{thm}

An analogous result holds for the multivariate case but the ``curvature matrix'' must be replaced by the self-adjoint algebra which the ``curvature matrices'' generate. In particular, there is an injective map from the lattice of reducing submodules to the lattice of reducing subspaces of this algebra. In general, this map is not surjective.

The curvature matrix having eigenvalues of multiplicity greater than one at a point doesn't imply the existence of other reducing submodules. However, such an assumption on an open set, under appropriate conditions, does. Analogous statements hold for the multivariate case.

\section{Powers of the Bergman Shift}\label{sec2.4}

Now let us consider the operator $M_{z^m}$ on the Bergman module $L^2_a({\bb D})$. Using $M_{z^m}$ as the contraction operator, we can define a quasi-free Hilbert module on ${\bb D}$ of multiplicity $m$, which we'll call ${\cl R}_m$. The question we want to consider is the determination of the reducing subspaces of $M_{z^m}$ or the reducing submodules of ${\cl R}_m$. Before we study the $m$th power of the Bergman shift, we consider the $m$th power of the unilateral shift; that is, $T_{z^m}$. Here the answer is straightforward since $T_{z^m}$ on $H^2({\bb D})$ is unitarily equivalent to $T_z\otimes I_m$ on $H^2({\bb D})\otimes {{\bb C}^m}$, since $T_{z^m}$ is an isometry of multiplicity $m$. In particular, there is a reducing subspace for every subspace of ${\bb C}^m$. Hence,   the lattice of reducing submodules is continuous and has infinitely many elements. The surprise, perhaps, is that this is not the case for $M_{z^m}$.

Let us begin with the case $m=2$. Then there are two obvious reducing submodules:\ ${\cl L}_0 = \bigvee\limits^\infty_{k=0} \{z^{2k}\}$, the span of the even powers of $z$, and ${\cl L}_1 = \bigvee\limits^\infty_{k=0} \{z^{2k+1}\}$, the span of the odd powers. Both ${\cl L}_0$ and ${\cl L}_1$ are invariant under multiplication by $z^2$ and together $\text{span } L^2_a({\bb D})$. Moreover, ${\cl L}_0$ and ${\cl L}_1$ are orthogonal since the $\{z^k\}$ form an orthogonal basis for $L^2_a({\bb D})$. Hence, $L^2_a({\bb D}) = {\cl L}_0 \oplus {\cl L}_1$ and ${\cl L}_0$ and ${\cl L}_1$ are reducing submodules. The question is whether there are other reducing submodules. As pointed out above, there are infinitely many more for the analogous construction using the Hardy module. However, for the Bergman module, these are all and the lattice of reducing submodules is discrete and finite. 

There are at least two ways to see that this is the case, one involving complex geometry and the other depends  directly on operator theory exploiting the fact that all these operators are unilateral weighted shifts. Let us consider the  complex geometric approach first.

Corresponding to the submodules ${\cl L}_0$ and ${\cl L}_1$, there are Hermitian holomorphic vector bundles with curvatures ${\cl K}_0(\omega) d\omega d\bar\omega$ and ${\cl K}_1(\omega)d\omega d\bar\omega$. If there were another anti-holomorphic sub-bundle of $E^*_{{\cl R}_2}$ whose ortho complement is also anti-holomorphic, then the curvature for $M_{z^2}$ on ${\cl R}_2$ would have to decompose as a direct sum corresponding to it as we discussed in the last section. In particular,  once one shows that the curvatures of $E^*_{{\cl L}_0}$ and $E^*_{{\cl L}_1}$ at 0 are distinct, such a decomposition is seen to be impossible and the result will follow.

Let $\gamma_\omega= \gamma^0_\omega$ be the vector in $L^2_a({\bb D})$ defined in \eqref{eq1.1} such that $M_z\gamma_\omega = \bar\omega \gamma_\omega$ for $\omega$ in ${\bb D}$. If $\eta$ and $-\eta$ are complex numbers such that $\eta^2=\omega$, then $\gamma_\eta$ and $\gamma_{-\eta}$ span the eigenspace for $M_{z^2}$ at $\omega$. If we differentiate the equation $M^*_z\gamma_\omega = \bar\omega\gamma_\omega$ using $\frac{d}{d\bar\omega}$, we obtain $M^*_z\gamma'_\omega = \gamma_\omega + \bar\omega \gamma'_\omega$, where $\gamma'_\omega = \frac{d}{d\bar\omega} \gamma_\omega$. An easy argument shows that a basis for $\text{KER}(M_{z^2}-\omega)^{*2}$ is given by $\gamma_\eta$, $\gamma'_\eta$, $\gamma_{-\eta}$ and $\gamma'_{-\eta}$ for $\omega\ne 0$. If one obtains the matrix for the nilpotent operator $(M_{z^2}-\omega)|_{\text{KER}(M_{z^2}-\omega)^{*2}}$, that would solve the problem. However, the calculation required is tedious. A better approach is to find anti-holomorphic sections for $E^*_{{\cl R}_2}$ which lie in ${\cl L}_0$ and ${\cl L}_1$. Another way to view this approach is to find global anti-holomorphic sections which also span the eigenspace at 0.

Defining $\nu_0$ as the sum of $\gamma_\eta$ and $\gamma_{-\eta}$ to obtain a section involving only the even powers of $z$ is more or less an obvious step. If one considers the difference of $\gamma_\eta$ and $\gamma_{-\eta}$, one obtains a section involving only the odd powers of $z$. However, the resulting section vanishes at 0 which is why one obtains $\nu_1$ by first dividing by $\eta$.

If we set $\nu_0 =\frac12(\gamma_\eta+\gamma_{-\eta})$ and $\nu_1 = \frac{\sqrt 2}{2\eta} (\gamma_\eta - \gamma_{-\eta})$, then easy calculations show that these sections form an anti-holomorphic orthogonal frame for $E^*_{{\cl R}_2}$.
Moreover, $\nu_0$ is a section for the bundle corresponding to ${\cl L}_0$ and $\nu_1$ is a section for ${\cl L}_1$. This follows since the Taylor series for $\nu_0$ involves only even powers of $z$ while the Taylor series for $\nu_1$ involves only the odd powers of $z$. The curvature is calculated using the formula $\bar\partial H^{-1}\partial H$, where $H$ is the Grammian matrix
\begin{equation}\label{eq2.6}
H(\pmb{\nu}) = \begin{pmatrix}
\langle\nu_0,\nu_0\rangle&\langle \nu_0,\nu_1\rangle\\
\langle \nu_1,\nu_0\rangle&\langle\nu_1,\nu_1\rangle
\end{pmatrix}\quad \text{and}\quad \pmb{\nu} = (\nu_0,\nu_1).
\end{equation}
Since $\langle \nu_0,\nu_1\rangle \equiv 0$, we obtain the diagonal matrix $\left(\begin{smallmatrix}\frac{1+|\omega|^2}{(1-|\omega|^2)^2}&0\\ 0&\frac1{(1-|\omega|^2)^2}\end{smallmatrix}\right)$.

A simple calculation yields the curvature matrix
\begin{equation}\label{eq2.7}
\begin{pmatrix}
-\frac{3+2|\omega|^2+3|\omega|^4}{(1-|\omega|^4)^2}&0\\
0&-\frac2{(1-|\omega|^2)^2}
\end{pmatrix}
\end{equation}
and hence at 0 we have $\left(\begin{smallmatrix} -3&0\\ 0&-2\end{smallmatrix}\right)$. Since the eigenvalues are distinct, the lattice of reducing submodules is discrete and finite. Moreover, since the curvature for the Bergman module is $-\frac2{(1-|\omega|^2)^2}$, we see that $M_{z^2}|_{{\cl L}_1}$ is unitarily equivalent to  $M_z$ on $L^2_a({\bb D})$.

This same approach, and  extended calculation, yields the same results for $M_{z^m}$ on $L^2_a({\bb D})$ for all $m$. That is, the lattice of reducing subspaces of $M_{z^m}$ is finite and discrete with the minimal elements $\{{\cl L}_{m,k}\}$, $0 \le k < m$ determined by the span of the powers of $z$ that are congruent to $k$ modulo $m$ and $M_{z^m}$ on ${\cl L}_{m,m-1}$ is unitarily equivalent to $M_z$ on $L^2_a({\bb D})$.

We won't provide the details for general $m$ but record the following result.

\begin{thm}\label{thm2.3}
For a positive integer $m>1$ and $0\le k<m$, set ${\cl L}_{m,k} = \vee \{z^\ell \colon \ \ell = k \bmod m\} \subset L^2_a({\bb D})$. Then $\{{\cl L}_{m,k}\}$ are reducing submodules for $M_{z^m}$ and the curvature for $M_{z^m}|_{{\cl L}_{m,k}}$ at 0 is $-\frac{m+k}k$. Thus these are all the reducing submodules for $M_{z^m}$ and the lattice of reducing submodules of ${\cl R}_m$ is finite and descrete.
\end{thm}

As a result, none of the operators $M_{z^m}|{\cl L}_{m,k}$ are unitarily equivalent. However, the above calculations can be used to show that ${\cl L}_{m,m-1}$ is unitarily equivalent to $L^2_a({\bb D})$. 

The approach to this question outlined above was obtained with J.M. Landsberg. The result was known and is part of a program begun by Zhu \cite{Zhu} and extended by Hu, Sun, Xu and Yu \cite{HSXY} and Sun, Zheng and Zhong \cite{SZZ} investigating the reducing subspaces of $M_B$ on $L^2_a({\bb D})$, for $B$ a finite Blaschke product.

As mentioned earlier, the operators obtained by restricting $M_{z^m}$ to the ${\cl L}_{m,k}$ are unilateral weighted shifts. It is a known result \cite{Shie} on how to determine when two such unilateral weighted shifts are unitarily equivalent or similar. Again let us provide the details for the $m=2$ case.

Since $\{z^k\}^\infty_{k=0}$ is an orthogonal basis, one obtains an orthonormal basis $\{e_k\}^\infty_{k=0}$, where $e_k = \frac{z^k}{\|z^k\|} = \sqrt{k+1} z^k$ and $M_ze_k = \sqrt{\frac{k+1}{k+2}} e_{k+1}$ for $k=0,1,2,\ldots$~. Further, the orthonormal bases for ${\cl L}_0$ and ${\cl L}_1$ are given by $\{e_{2k+1}\}^\infty_{k=0}$ and $\{e_{2k}\}^\infty_{k=0}$, respectively. Thus, if $T_0 = M_{z^2}|_{{\cl L}_0}$ and $T_1 = M_{z^2}|_{{\cl L}_1}$, then
\begin{align*}
T_0e_{2k} &= T_0(\|z^{2k}\|z^{2k}) = \frac{\|z^{2k}\|}{\|z^{2k+2}\|} \|z^{2k+2}\| z^{2k+2} = \sqrt{\frac{2k+1}{2k+3}} e_{2k+2}\\
\intertext{and}
T_1e_{2k+1} &= T_1(\|z^{2k+1}\|z^{2k+1}) = \frac{\|z^{2k+1}\|}{\|z^{2k+3}\|} \|z^{2k+3}\|z^{2k+3} = \sqrt{\frac{2k+2}{2k+4}} e_{2k+3}\\
&= \sqrt{\frac{k+1}{k+2}} e_{2k+3}. 
\end{align*}
Thus since the weights are not equal, we see that $T_0$ on ${\cl L}_0$ and $T_1$ on ${\cl L}_1$ are not unitarily equivalent but since the weights do agree, we see that $T_1$ on ${\cl L}_1$ is unitarily equivalent to $M_{z}$ on $L^2_a({\bb D})$. More specifically, if we define $X\colon \ L^2_a({\bb D})\to {\cl L}_1$ so that $Xe_k = e_{2k+1}$ for $k=0,1,2,\ldots$, then $X$ is a unitary operator satisfying $XM_{z} = T_1X$.

If we define $Y\colon \ L^2_a({\bb D})\to {\cl L}_0$ so that $Ye_k = c_k e_{2k}$, where $c_k = \sqrt{\frac{2(2k)(2k+1)}{3(k+1)(k+2)}}$, then $YM_z = T_0Y$. Moreover, $Y$ is bounded and invertible since $\lim\limits_{k\to\infty} c_k = \frac83 >0$. Thus the operators obtained by restricting $M_{z^2}$  to both ${\cl L}_0$ and ${\cl L}_1$ are similar to $M_z$ on $L^2_a({\bb D})$.

The same kind of argument establishes the analogous result for arbitrary $m>1$.

Finally, let us observe that the similarity  results show that $M_{z^m}$ on $L^2_a({\bb D})$ is similar to $M_z\otimes I_m$ on $L^2_a({\bb D})\otimes {\bb C}^m$. This result raises an interesting question related to the program of Zhu mentioned earlier.\medskip 

\n {\bf Question 6.} Is $M_B$ on $L^2_a({\bb D})$  similar to $M_z\otimes I_m$ on $L^2_a({\bb D}) \otimes {\bb C}^m$, where $m$ is the multiplicity of the finite Blaschke product $B(z)$?

\section{Submodules Obtained from Ideals}\label{sec2.5}

As mentioned earlier, the problem of describing submodules of a quasi-free Hilbert module is a nearly impossible one in most cases although much has been learned in the past decade or two. For example, consider what is known about the submodules for the Hardy module \cite{Beur} and the Bergman modules (cf. \cite{HKZ}). For the multivariate case, we will illustrate the wide variety of possibilities with another rigidity theorem and an example.

Let ${\cl R}$ be a quasi-free Hilbert module over ${\bb B}^n$ of multiplicity one. If for $I$ an ideal in ${\bb C}[\pmb{z}]$, $[I]_{\cl R}$ denotes the closure of $I$ in ${\cl R}$, then $[I]_{\cl R}$ is a submodule of ${\cl R}$. Of course, not all submodules arise in this manner as can be seen by considering submodules of the Hardy module $H^2({\bb D})$. However, the question in which we are interested is when is $[I]_{\cl R} \cong [\widetilde I]_{\cl R}$ for two ideals $I$ and $\widetilde I$ of ${\bb C}[\pmb{z}]$. There is now  considerable literature on this question \cite{D-P}, \cite{DPSY}, \cite{C-G}. We state one of the early fundamental results.

\begin{thm}\label{thm2.4}
Suppose ${\cl R}$ is a quasi-free Hilbert module over ${\bb B}^n$ of multiplicity one and $I, \widetilde I$ are ideals in ${\bb C}[\pmb{z}]$. Assume further that
\begin{equation}\label{eq2.8}
\begin{array}{ll}
\text{\rm 1)}& \text{the codimension of each of the algebraic components }\\
&\text{of ${\cl Z}(I)$ and ${\cl Z}(\widetilde I)$ is greater than one and}\\
\text{\rm 2)}&\text{these same components all intersect ${\bb B}^n$.}
\end{array}
\end{equation}
Then $[I]_{\cl R}$ and $[\widetilde I]_{\cl R}$ are quasi-similar iff $I = \widetilde I$.
\end{thm}

Here, ${\cl Z}(I)$ denotes the zero variety of $I$ or the set of common zeros of the polynomials in $I$ and by algebraic component of $I$ is meant the zero variety of a prime ideal in a primary decomposition of $I$. If we omit (2), then $[I]_{\cl R}$ and $[\widetilde I]_{\cl R}$ quasi-similar holds iff $[I]_{\cl R} = [\widetilde I]_{\cl R}$.

We conclude with one example which illustrates the kind of techniques used which trace back to results of Zariski and Grothendieck.

We consider the Hardy module $H^2({\bb D}^2)$ over the bidisk which is the closure of ${\bb C}[z_1,z_2]$ in $L^2({\bb T}^2)$ and the submodule $H^2_0({\bb D}^2)$ of functions in $H^2({\bb D}^2)$ that vanish at the origin $\pmb{0}$. A simple calculation shows that
\begin{equation}\label{eq2.9}
\begin{split}
H^2({\bb D}^2)/I_{\bf 0}\cdot H^2({\bb D}^2) &\cong {\bb C}_{\bf 0}, \quad \text{while}\\
H^2_0({\bb D}^2)/I_{\bf 0}\cdot H^2_0({\bb D}^2) &\cong {\bb C}_{\bf 0} \oplus {\bb C}_0.
\end{split}
\end{equation}
Here ${\bb C}_{\bf 0}$ denotes the Hilbert module over ${\bb C}[\pmb{z}]$ for the Hilbert space ${\bb C}$ in which module multiplication is defined $p\cdot\lambda = p({\bf 0})\lambda$ for $p$ in ${\bb C}[\pmb{z}]$ and $\lambda$ in ${\bb C}$.

The same results hold for $H^2({\bb B}^2)$ but the calculations for the bidisk are more transparent.

If $X$ and $Y$ define a quasi-similarity between $H^2({\bb D}^2)$ and $H^2_{\bf 0}({\bb D}^2)$, then the localized maps $\widehat X(\pmb{0})$ and $\widehat Y(\pmb{0})$ would define an isomorphism between ${\bb C}_{\bf 0}$ and ${\bb C}_{\bf 0} \oplus {\bb C}_{\bf 0}$ which is impossible.

Note that one can show if ${\cl R}$ and $\widetilde{\cl R}$ are quasi-free Hilbert modules of multiplicity one over $A({\bb B}^n)$ and $I$ and $\widetilde I$ are ideals in ${\bb C}[\pmb{z}]$ satisfying \eqref{eq2.8} such that $[I]_{\cl R}$ and $[\widetilde I]_{\widetilde{\cl R}}$ are quasi-similar or, $[I]_{\cl R}\sim [\widetilde I]_{\cl R}$, then $I=\widetilde I$. However, it is not true that $[I]_{\cl R}\sim [I]_{\widetilde{\cl R}}$ implies ${\cl R}\sim \widetilde{\cl R}$, but it seems reasonable to ask when it does.  For an example  when this implication does not hold, consider the Hilbert module $L^2_a(\mu)$ over $A({\bb D})$ obtained from the completion of ${\bb C}[z]$ in $L^2(\mu)$, where $\mu$ is Lebesgue measure on ${\bb D}$ plus the point mass at 0 and $I$ is the principal ideal generated by $z$. Then the closures of $I$ in $L^2_a(\mu)$ and $L^2_a({\bb D})$ are unitarily equivalent but $L^2_a(\mu)$ and $L^2_a({\bb D})$ are not even quasi-similar.\medskip  

\n {\bf Question 7.} Find conditions on an ideal $I$ in ${\bb C}[\pmb{z}]$ and/or on the quasi-free Hilbert modules ${\cl R}$ and $\widetilde{\cl R}$ of multiplicity one over $A({\bb B}^n)$ so that $[I]_{\cl R}$ quasi-similar to $[I]_{\widetilde{\cl R}}$ implies that ${\cl R}$ and $\widetilde{\cl R}$ are quasi-similar. What about the same question for similarity or unitary equivalence?

\clearpage\thispagestyle{empty}\cleardoublepage

\chapter{THIRD LECTURE}\label{chp3}

\section{Isomorphic Submodules}\label{sec3.1}

In this last lecture we take up a new topic which we believe to have intrinsic interest and to provide an impetus to develop new tools for the study of multivariate operator theory.

If ${\cl L}$ is a non-zero submodule of $H^2({\bb D})$, then a consequence of Beurling's Theorem \cite{Beur} is that ${\cl L}$ is isometrically isomorphic to $H^2({\bb D})$ itself. In particular, there exists an inner function $\theta$ so that ${\cl L} = \theta H^2({\bb D})$ and the operator $T_{\theta}$ defines an isometry on $H^2({\bb D})$ which establishes a module isometric isomorphism between $H^2({\bb D})$ and ${\cl L}$. On the other hand, a result of Richter \cite{Rich} shows that if ${\cl L}$ is a submodule of $L^2_a({\bb D})$ which is isometrically isomorphic to $L^2_a({\bb D})$, then ${\cl L} = L^2_a({\bb D})$. Thus, in one case every non-zero submodule is isometrically isomorphic to the module itself while in the other, no proper submodule is. J.~Sarkar and the author recently investigated this phenomenon in \cite{DS} and it is some of the results from \cite{DS} we discuss here. In this lecture, we will allow Hilbert modules over $A(\Omega)$ for bounded domains $\Omega\subset {\bb C}^n$.

Our first result shows that the existence of a proper submodule of \emph{finite} codimension, isometrically isomorphic to the original module is a single operator or one variable phenomenon; that is, it occurs only for Hilbert modules over $A(\Omega)$ for $\Omega \subset {\bb C}$. Before proceeding, however, we need conditions ensuring the ``fullness'' of the submodule.

If ${\cl R}$ is a Hilbert module over $A(\Omega)$, then there exists a submodule ${\cl L}$ of ${\cl R}$ isometrically isomorphic to ${\cl R}$ iff there exists an isometric module map $V$ on ${\cl R}$ with range ${\cl L}$. If $V$ is unitary, then ${\cl L} = {\cl R}$. We are interested in the opposite extreme; that is, when $V$ is pure or $\bigcap\limits^\infty_{k=0} V^k{\cl R} = (0)$ in which case we have the following result.

\begin{thm}\label{thm3.1}
Let ${\cl R}$ be a quasi-free Hilbert module of finite multiplicity over $A(\Omega)$ for $\Omega$ a bounded domain in ${\bb C}^n$. If there exists a submodule ${\cl L}$ of ${\cl R}$ such that $\text{\rm DIM}_{\bb C} {\cl R}/{\cl L} <\infty$ and the corresponding isometric module map is pure, then $n=1$ and $\Omega\subset {\bb C}$. 
\end{thm}

Our approach to prove this result  involves the Hilbert--Samuel polynomial.
Recall that the Hilbert--Samuel polynomial for the Hilbert module ${\cl M}$ is a polynomial $h^{\cl M}_{\pmb{\omega}}(z)$ in one variable so that $h^{\cl M}_{\pmb{\omega}}(k) = \text{DIM}_{\bb C} ~{\cl M}/[I^k_{\pmb{\omega}}\cdot {\cl M}]$ for $k\gg 0$. A necessary condition for the existence of such a polynomial  is that $\text{DIM}_{\bb C} {\cl M}/[I_{\pmb{\omega}}\cdot {\cl M}]<\infty$ and in \cite{D-Y} K.\ Yan and the author established  the existence of such a polynomial under this assumption extending the earlier work of Hilbert and Samuel to this context. More recently, Arveson has also considered this notion \cite{A1} in his study of Hilbert modules related to $H^2_n$.

One can show by analyzing local frames for the Hermitian holomorphic bundle $E_{\cl R}$ for ${\cl R}$ a quasi-free Hilbert module of finite multiplicity $m$, that the Hilbert--Samuel polynomial for ${\cl R}$ does not depend on which quasi-free Hilbert module is chosen. We consider the argument to establish this independence for the case $n=2$.

Assume ${\cl R}$ is a quasi-free Hilbert  module over $A(\Omega)$ such that $\text{DIM } {\cl R}/[I_{\pmb{\omega}}\cdot$ ${\cl R}] = m<\infty$ and $\Omega\subset {\bb C}^2$. For some neighborhood $V$ of $\pmb{\omega}_0$, there exist vectors $\{f_i\}^m_{i=1}$ such that $\{\hat f_i(\pmb{\omega})\}^m_{i=1}$  is a basis for ${\cl R}_{\pmb{\omega}}$ for $\pmb{\omega}$ in $V$. Then one can show that $\{\frac{\partial^{\alpha_1+\alpha_2}}{\partial z^{\alpha_1}_1 \partial z^{\alpha_2}_2} \hat f_i(\pmb{\omega})\}_{0\le\alpha_1+\alpha_2<k}$; $1\le i\le m$ is a basis for ${\cl R}/[{I^k_{\pmb{\omega}}\cdot{\cl R}}]$ for $\pmb{\omega}$ in $V$. However, the cardinality of the set $\{(\alpha_1,\alpha_2)\colon \ \alpha_1,\alpha_2\ge 0$ and $\alpha_1 + \alpha_2<k\}$ is $\frac12(k-1)^2$. Therefore, $h^{\cl R}_{\pmb{\omega}}$ has degree 2. Hence, in a similar manner one can show for $\Omega\subset {\bb C}^n$ that $h^{\cl R}_{\pmb{\omega}}$ has degree $n$. The proof of the theorem is completed by showing that under the hypotheses of an isometrically isomorphic submodule of finite codimension,  the Hilbert--Samuel polynomial is actually linear, or that $n=1$. We first need a lemma which  is a key step in the analysis of Hilbert modules containing a pure isometrically isometric submodule.

\begin{lem}\label{lem3.1}
Let ${\cl M}$ be a Hilbert module over ${\bb C}[\pmb{z}]$ with a pure isometrically isomorphic submodule ${\cl L}$. Then there exists a Hilbert space ${\cl E}$; $\varphi_1,\ldots, \varphi_n$ in $H^\infty_{{\cl L}({\cl E})}({\bb D})$ such that $[\varphi_i(\omega), \varphi_j(\omega)] = 0$ for $\omega$ in ${\bb D}$ and $1\le i,j\le n$; and an isometrical isomorphism $\Psi\colon \ {\cl M}\to H^2_{\cl E}({\bb D})$  such that $\Psi$ is a module isomorphism where $M_p$ on $H^2_{\cl E}({\bb D})$ is defined to be multiplication by $p(\varphi_1,\ldots,\varphi_n)$. Moreover, $T_z$ is the isometric module map that defines ${\cl L}$.
\end{lem}

As is perhaps obvious, to establish the lemma one uses the von Neumann--Wold decomposition to identify $V$ on ${\cl M}$ with the Toeplitz operator $T_z$ on $H^2_{\cl E}({\bb D})$, where ${\cl E} = {\cl M}/V{\cl M}$, and then uses the fact that all operators that commute with $T_z$ have the form $T_\eta$ for $\eta$ in $H^\infty_{{\cl L}({\cl E})}({\bb D})$.

\begin{lem}\label{lem3.2}
Let ${\cl M}$ be a Hilbert module over ${\bb C}[\pmb{z}]$ with a pure isometrically isomorphic submodule ${\cl L}$ such that $\text{\rm DIM}_{\bb C} {\cl M}/{\cl L} < \infty$ and $\text{\rm DIM}_{\bb C} {\cl M}/[{I_{\pmb{\omega}_0}\cdot {\cl M}}] < \infty$ for some $\pmb{\omega}_0$. Then $h^m_{\pmb{\omega}_0}$ is linear.
\end{lem}

\begin{proof}
Without loss of generality we can assume that $\pmb{\omega}_0 = \pmb{0}$. Applying the previous lemma, we can identify ${\cl M}$ with $H^2_{\cl E}({\bb D})$ using $\Psi$, where ${\cl E} = {\cl M}/{\cl L}$ is finite dimensional, with the module action on $H^2_{\cl E}({\bb D})$  defined by an $n$-tuple $\{\varphi_i\}^n_{i=1}$ of functions in $H^\infty_{{\cl L}({\cl E})}({\bb D})$ such that $[\varphi_i(z), \varphi_j(z)] = 0$ for $z$ in ${\bb D}$ and $1\le i,j\le n$ and ${\cl L} = T_zH^2_{\cl E}({\bb D})$. Since $\text{DIM}_{\bb C} {\cl M}/[I_{\pmb{0}}\cdot {\cl M}] < \infty$, it follows that the closure of the $T_{\varphi_1}\cdot H^2_{\cl E}({\bb D}) +\cdots+ T_{\varphi_n}\cdot H^2_{\cl E}({\bb D})$ is a closed subspace of $H^2_{\cl E}({\bb D})$ having finite codimension and which is invariant under $T_z$. Thus by the Beurling--Lax--Halmos Theorem (cf.\ \cite{Sz-NF}) it follows that there exists an inner function $\Theta(z)$ in $H^\infty_{{\cl L}({\cl E})}({\bb D})$ so that $\Theta H^2_{\cl E}({\bb D}) = \text{clos}[T_{\varphi_1}\cdot H^2_{\cl E}({\bb D}) +\cdots+ T_{\varphi_n}\cdot H^2_{\cl E}({\bb D})]$. Moreover, since $\Theta H^2_{\cl E}({\bb D})$ has finite codimension in $H^2_{\cl E}({\bb D})$, it follows that $\Theta(z)$ is a rational function and $\Theta(e^{it})$ is unitary for $e^{it}$ in ${\bb T}$ (cf.\ \cite{Sz-NF}). Moreover, one can define $\det \Theta$, the determinant of $\Theta(z)$, on ${\bb D}$ to obtain a rational scalar-valued inner function satisfying 
\[
(\det\Theta) H^2_{\cl E}({\bb D}) \subseteq \Theta H^2_{\cl E}({\bb D}).
\]
The proof of this fact depends on Cramer's Rule and is due to Helson \cite{Hel}. Now we have
\[
[I^2_{\pmb{0}}\cdot {\cl M}] = \bigvee_{|\pmb{\alpha}|=2} \pmb{z}^{\pmb{\alpha}}\cdot {\cl M}, \text{ where } |\pmb{\alpha}| = \alpha_1+\alpha_2.
\]
Hence
\begin{align*}
\bigvee_{1\le i,j\le n} T_{\varphi_i}T_{\varphi_j} H^2_{\cl E}({\bb D}) &= T_{\varphi_1} \Theta H^2_{\cl E}({\bb D}) +\cdots+ T_{\varphi_n} \Theta H^2_{\cl E}({\bb D})\\
&\supseteq T_{\varphi_1} \det \Theta H^2_{\cl E}({\bb D}) +\cdots+ T_{\varphi_n} \det \Theta H^2_{\cl E}({\bb D})\\
&\supseteq \det \Theta(T_{\varphi_1} H^2_{\cl E}({\bb D}) +\cdots+ T_{\varphi_n} H^2_{\cl E}({\bb D}))\\
&\supseteq \det \Theta \Theta H^2_{\cl E}({\bb D}) \supseteq (\det\Theta)^2 H^2_{\cl E}({\bb D}).
\end{align*}
By induction, one obtains
\[
(\det\Theta)^k H^2_{\cl E}({\bb D}) \subseteq \bigvee_{1\le i_1,i_2,\ldots,i_k\le n} T_{\varphi_{i_1}} T_{\varphi_{i_2}}\cdots T_{\varphi_{i_k}} H^2_{\cl E}({\bb D}).
\]
Thus using $\Psi$ we have
\[
\text{DIM}_{\bb C} {\cl M}/[I^k_{\pmb{0}}{\cl M}] \le \text{DIM}_{\bb C} H^2_{\cl E}({\bb D})/(\det\Theta)^k H^2_{\cl E}({\bb D}).
\]
Therefore, $h^m_{\pmb{\omega}}(k) \le k\ell\cdot \text{DIM}_{\bb C} {\cl E}$, where $\ell$ is the number of zeros of $\det\Theta$ counted multiply, and hence $h^m_{\pmb{\omega}}$ is linear.
\end{proof}

Combining the lemmas, we complete the proof of the theorem. Hence, to study isometrically isomorphic submodules of finite codimension, we can assume that $\Omega$ is a domain in ${\bb C}$. If we assume further that $\Omega$ is ``nice,'' then we can characterize the situation completely.

\begin{thm}\label{thm3.2}
If ${\cl R}$ is a contractive, quasi-free Hilbert module over $A({\bb D})$ of finite multiplicity containing a pure isometrically isomorphic submodule of finite codimension, then ${\cl R}$ is isometrically isomorphic to $H^2_{\cl E}({\bb D})$ as $A({\bb D})$-Hilbert modules with $\text{\rm DIM}_{\bb C} {\cl E}$ finite and equal to the multiplicity of ${\cl R}$.
\end{thm}

\begin{proof}
Using Lemma \ref{lem3.1}, one reduces the question to a single Toeplitz operator $T_\varphi$ on some $H^2_{\cl E}({\bb D})$, where $\varphi$ is in $H^\infty_{{\cl L }({\cl E})}({\bb D})$ and $\text{DIM}_{\bb C} {\cl E} < \infty$. But one knows since ${\cl R}$ is contractive, that $\|\varphi(z)\| \le 1$ for $z$ in ${\bb D}$. Moreover, since ${\cl R}$ is quasi-free, we can conclude that the spectrum of the matrix $\varphi(e^{it})$ is contained in ${\bb T}$ a.e. Since we are on the finite dimensional space ${\cl E}$, it follows that $\varphi(e^{it})$ is unitary a.e. Hence $T_\varphi$ is an isometry and the module action it defines on $H^2_{\cl E}({\bb D})$ yields a module which is isomorphic to a Hardy module $H^2_{\cl F}({\bb D})$. Finally, the multiplicity of $T_\varphi$, and hence the dimension of ${\cl F}$, is the dimension of ${\cl R}/I_0\cdot {\cl R}$ which is the multiplicity of ${\cl R}$ as a quasi-free Hilbert module which was assumed to be finite.

A key step in the argument is to observe that $T_\varphi-\omega$ Fredholm implies that $\varphi(e^{it})-\omega$ is bounded away from 0 for $\omega$ in ${\bb D}$. 
\end{proof}

\n {\bf Question 8.} Does the theorem hold with the assumption that ${\cl R}$ is only weakly quasi-free or without assuming that $I_{\pmb{\omega}}\cdot {\cl R}$ is closed for $\omega$ in ${\bb D}$?

If we assume that $\Omega$ is a finitely connected domain in ${\bb C}$ for which $\partial\Omega$ consists of a finite number of simple closed curves, then a similar argument yields the analogous result with Hilbert modules defined by bundle shifts replacing the Hardy modules on ${\bb D}$.

Recall that if $\pi^1(\Omega)$ is the fundamental group of $\Omega$ and $\alpha\colon \ \pi^1(\Omega)\to {\cl U}({\cl E})$ is a unitary representation, then the bundle shift $H^2_\alpha(\Omega)$ is defined as a space of holomorphic sections of the flat Hermitian holomorphic bundle determined by $\alpha$. The multiplicity of $H^2_\alpha(\Omega)$ equals the dimension of ${\cl E}$. The norm can be defined using harmonic measure on $\Omega$ for some point in $\Omega$. The theory of such Hilbert modules was developed by Abrahamse and the author \cite{A-D} and they play the role for multiply connected domains that the Hardy module does for the unit disk. If one specializes to the case of the annulus, one obtains the spaces of modulus automorphic holomorphic functions studied by Sarason \cite{Sara}.

\begin{thm}\label{thm3.4}
Let $\Omega$ be a finitely connected domain with $\partial\Omega$ consisting of simple closed curves. If ${\cl R}$ is a quasi-free Hilbert module over $A(\Omega)$ of finite multiplicity for which there exists a pure isometrically isomorphic submodule of finite codimension, then ${\cl R}$ is isomorphic to $H^2_\alpha(\Omega)$ for some unitary representation $\alpha\colon \ \pi^1(\Omega)\to {\cl U}({\cl E})$, where $\text{DIM } {\cl E}$ is finite and equals the multiplicity of ${\cl R}$.
\end{thm}

Bundle shifts are examples of ``subnormal modules.'' Recall that a Hilbert module ${\cl S}$ is \emph{subnormal} if there exists a \emph{reductive} Hilbert module ${\cl N}$ or one for which module multiplication defines commuting normal operators, and such that ${\cl S}$ is a submodule of ${\cl N}$. If ${\cl N}$ extends to a module over $C(\partial A)$, where $\partial A$ denotes the \v Silov boundary of the algebra $A(\Omega)$ over which ${\cl L}$ is a Hilbert module, then ${\cl L}$ is said to be a \emph{\v Silov module}. Bundle shifts are examples of \v Silov modules as are the Hardy modules.\medskip 

\n {\bf Question 9.} If ${\cl R}$ is a weakly quasi-free Hilbert module of finite multiplicity containing a pure isometrically isomorphic  submodule, must it be a \v Silov module? Conversely, does every \v Silov module possess such a submodule, at least for $\Omega\subset {\bb C}$?

\begin{rem}\label{rem3.1}
It seems likely that a \v Silov module which is weakly quasi-free of finite multiplicity is actually  quasi-free. If so, then one could add that hypothesis to the question.
\end{rem}

There is another question which one can ask in this context but which is not related to the topic of this section.\medskip 

\n {\bf Question 10.} Let $H^2_\alpha(\Omega)$ be a bundle shift of finite multiplicity over $\Omega$. What is the precise relation between the flat unitary holomorphic bundle $E_\alpha$ determined by $\alpha$ and the bundle $E_{H^2_\alpha(\Omega)}$? Does $E_\alpha$ represent the holonomy of $E_{H^2_\alpha(\Omega)}$ in some sense?

\section{The Case of Infinite Codimension}\label{sec3.2}

What can one say about a (weakly) quasi-free Hilbert module of finite multiplicity that contains a pure isometrically isomorphic submodule of infinite codimension? We offer a couple of observations before proceeding to some results.

If ${\cl R}$ is a (weakly) quasi-free Hilbert module over $A(\Omega)$, then ${\cl R} \otimes H^2({\bb D})$ is a (weakly) quasi-free Hilbert module over $A(\Omega\times {\bb D})$ and ${\cl R} \otimes H^2_0({\bb D})$ is a pure isometrically isomorphic submodule. Note that the boundary of $\Omega\times {\bb D}$ has ``corners'' or is not smooth.\medskip 

\n {\bf Question 11.} Could we say something about ${\cl R}$ if we assume that $\Omega$ is strongly pseudo-convex with smooth boundary?

The nicest domain in ${\bb C}^n$, of course, is the unit ball ${\bb B}^n$. The existence of an inner function, established by Aleksandrov \cite{Alek}, yields a pure isometrically isomorphic submodule $\theta H^2({\bb B}^n)$ of $H^2({\bb B}^n)$.

If we add the assumption that the module is essentially reductive, we can reach some rather surprising conclusions. Recall that a Hilbert module ${\cl M}$ is  said to be \emph{essentially reductive} if all of the operators on ${\cl M}$ defined by  module multiplication are essentially normal.

\begin{thm}\label{thm3.5}
If ${\cl R}$ is an essentially reductive Hilbert module over $A(\Omega)$ for some bounded domain $\Omega$ of ${\bb C}^n$ that contains a pure isometrically isomorphic submodule, then ${\cl R}$ is subnormal.
\end{thm}

\begin{proof}
Again we use the representation of ${\cl R}$ as $H^2_{\cl E}({\bb D})$ with module multiplication defined by a commuting $n$-tuple  $\{\varphi_i\}^n_{i=1}$ of commuting functions in $H^\infty_{{\cl L}({\cl E})}({\bb D})$. Of course, we must allow ${\cl E}$ to have infinite dimension in this case. However, the assumption that the algebra generated by the operators $T_{\varphi_i}$ is essentially normal implies that the operators $\varphi_i(e^{it})^*$ and $\varphi_j(e^{it})$ commute for $e^{it}$ in ${\bb T}$ a.e.\ for $1\le i,j\le n$. The key step in this argument is to observe that $[T^*_{\varphi_i},T_{\varphi_j}]$ compact implies that $[L^*_{\varphi_i}, L_{\varphi_j}]$ is compact.  Then one notes that the $C^*$-algebra generated by $\{L_{\varphi_i}\}$ contains no non-zero compact operator. Thus $\{\varphi_i(e^{it})\}^n_{i=1}$ is an $n$-tuple of commuting normal operators for $e^{it}$ in ${\bb T}$ a.e. Hence, the Hilbert module for $L^2_{\cl E}({\bb T})$ with module multiplication defined by $\{L_{\varphi_i}\}^n_{i=1}$ is reductive and extends the Hilbert module $H^2_{\cl E}({\bb D})$ with module multiplication defined by $\{T_{\varphi_i}\}^n_{i=1}$. Thus, the latter Hilbert module is subnormal and hence so is ${\cl R}$.
\end{proof}

In \cite{DS} we reprove a result of Chen and Guo \cite{C-G} that no proper submodule of $H^2_n$, for $n>1$, is isometrically isomorphic to $H^2_n$. The proof requires handling the possibility that the isometry, whose range is the submodule, is not pure.
Note that this result along with the  theorem provides another proof of the result of Arveson \cite{A2} that coordinate multiplication operators on the $n$-shift space are not jointly subnormal. Actually, this result was established earlier by Lubin in \cite{L} where he defined the space $H^2_m$ as a commuting weighted shift space. His purpose was to exhibit commuting subnormal operators, namely multiplication by the coordinate functions $M_{z_1},\ldots, M_{z_n}$, for which their sum and products are not subnormal. Hence there is no common normal extension.

One should note also based on formula \eqref{eq1.1} that the restriction of $M_{z_1}$ to cyclic subspace generated by $z^\ell_2$ is unitarily equivalent to a weighted Bergman shift for each $\ell$. Hence, $M_{z_1}$, on $H^2_2$, is the orthogonal direct sum of subnormal operators and hence is subnormal. The same is true for $M_{z_2}$. Therefore, $M_{z_1}$ and $M_{z_2}$ are commuting subnormal operators which have no joint normal extension. 

Further, we can use this theorem and the result of Athavale \cite{A} to prove an analogue of Theorem \ref{thm2.4} for the unit ball.

\begin{thm}\label{thm3.6}
If ${\cl R}$ is an essentially reductive  quasi-free Hilbert module over $A({\bb B}^n)$ which contains a pure isometrically isomorphic submodule ${\cl M}$, then ${\cl R}$ is isomorphic to $H^2_{\cl E}({\bb B}^n)$ and ${\cl M} = \theta H^2_{\cl E}({\bb B}^n)$ for some inner function $\theta$ and $\text{\rm DIM}_{\bb C} {\cl E} < \infty$.
\end{thm}

The proof is  analogous to  the earlier one using facts about row contractions on a Hilbert space of finite dimension. Note again that $H^2_{\cl E}({\bb B}^n)$ is a \v Silov module since $H^2_{\cl E}({\bb B}^n) \subseteq L^2_{\cl E}(\partial{\bb B}^n)$ and $\partial {\bb B}^n$ is the \v Silov boundary of $A({\bb B}^n)$.

In \cite{DS} we obtain more results about which subnormal Hilbert modules contain pure isometrically isomorphic submodules. In particular, we extend the results of Richter \cite{Rich} and Putinar \cite{Put} to Bergman modules over other domains.

\begin{thm}\label{thm3.7}
Let $\Omega$ be a bounded domain in ${\bb C}^n,\mu$ be a probability measure on $\Omega$ and $L^2_a(\mu)$ the closure of $A(\Omega)$ in $L^2(\mu)$. No proper submodule of $L^2_a(\mu)$ is isometrically isomorphic to $L^2_a(\mu)$.
\end{thm}

\vspace{.5in}

\n Department of Mathematics\\
Texas A\&M University\\
rdouglas@math.tamu.edu

\end{document}